\documentclass[reqno,11pt]{article}

\usepackage{amsmath,latexsym,amssymb}

\def\sqr#1#2{{\vcenter{\hrule height.#2pt

        \hbox{\vrule width.#2pt height#1pt \kern#1pt

                \vrule width.#2pt}

        \hrule height.#2pt}}}

\def\square{\mathchoice\sqr64\sqr64\sqr{4}3\sqr{3}3}

\def\QED{\hfill$\square$\break}

\def\demo{\noindent{\bf Proof: }}

\def\NN{{\mathbb  N }} 
\def\opn#1#2{\def#1{\operatorname{#2}}} 
\opn\Sq{Sq}
\opn\GCD{GCD}
\opn\Ker{Ker} 
\opn\Tor{Tor} 
\opn\reg{reg} 
\opn\chara{char}

\newtheorem{Theorem}{\sc Theorem}[section]
\newtheorem{Lemma}[Theorem]{\sc Lemma}

\newtheorem{Remark}[Theorem]{\sc Remark}
\newtheorem{Example}[Theorem]{\sc Example}
\newtheorem{Definition}[Theorem]{\sc Definition}

\newtheorem{Hint}[Theorem]{\sc Hint}

\begin{document}

\baselineskip=13pt

\pagestyle{empty}

\ \vspace{1.7in}

\noindent {\LARGE\bf Regularity jumps for powers of ideals
}

\vspace{.25in}

\noindent Aldo  \ Conca,   \\  Dipartimento di Matematica,  Universita' di
Genova, \\ Via Dodecaneso 35,  I-16146 Genova, Italia.   \\ {\it E-mail}: {\tt
conca@dima.unige.it}

\vspace{2.4cm}

\section{Introduction \hfill\break}

The Castelnuovo-Mumford regularity $\reg(I)$ is  one of the most
important invariants of a homogeneous  ideal  $I$ in a polynomial ring.   
A basic
question is how the regularity behaves with respect to taking powers of ideals. It is
known that in the long-run $\reg(I^k)$ is a linear function of $k$. We show that in the
short-run the regularity of
$I^k$  can be quite ``irregular".  For any given integer $d>1$ we
construct an ideal $J$  generated by $d+5$ monomials of degree $d+1$  in $4$ variables
such that  $\reg(J^k)=k(d+1)$ for every
$k<d$ and $\reg(J^d)\geq d(d+1)+d-1$.

\section{Generalities \hfill\break } \label{ng}

Let $K$ be a field. Let $R=K[x_1,\ldots,x_n]$ be the polynomial ring over
$K$.   Let
$I=\oplus_{i\in \NN}  I_i$ be a  homogeneous ideal. For every $i,j\in
\NN$   one defines the $ij$th graded Betti number of $I$ as
$$\beta_{ij}(I)=\dim_K \Tor^R_i(I,K)_j$$ and set 
$$t_i(I)=\max\{ j \,  | \,    \beta_{ij}(I)\neq 0\}$$  
 with  $t_i(I)=-\infty$ if  it happens that
$\Tor^R_i(I,K)=0$. The Castelnuovo-Mumford regularity 
$\reg(I)$ of $I$ is  defined as 
$$\reg(I)=\sup \{ t_i(I)-i \,  : i\in \NN \}. $$  By construction
$t_0(I)$ is the largest degree of a minimal generator of $I$.  The 
initial  degree  of  $I$ is the smallest  degree of a minimal generator
of $I$, i.e. it is the least index $i$ such that $I_i\neq 0$. The ideal 
$I$ has a   linear  resolution  if  its regularity is equal to its
initial degree. In other words,
$I$ has a  linear resolution if  its minimal generators all have  the same
degree and    the non-zero entries of the matrices  of the minimal free resolution of
$I$ all have   degree $1$.

The highest degree of a  generator  the $k$-th power  $I^k$  of $I$ is
bounded  above by
$k$ times the highest degree of a generator of $I$, i.e. 
$t_0(I^k)\leq k t_0(I)$.  One may wonder whether  the same
relation holds also for the Castelnuovo-Mumford regularity, that is,
whether the inequality 
$$  \reg(I^k)\leq k\reg(I) \eqno(1) $$  holds for every $k$. For some classes of ideals
$(1)$ holds, see e.g. \cite{CH} but in general it   does not. 
In Section 3 we present some examples of  ideals with linear resolution
whose square does not have a linear resolution.  On the other hand, it is known that for
every ideal
$I$ one has  $$\reg(I^k)=a(I)k+b(I) \quad \mbox{ for }\quad   \mbox{ for all } k\geq
c(I)   
\eqno(2)$$  where $a(I), b(I)$ and $c(I)$  are integers.   Cutkosky, Herzog
and  Trung \cite{CHT},  and  Kodiyalam \cite{K} proved independently that  $(2)$ holds.
They also shown that  $a(I)$ is bounded above by  the largest degree of a generator of
$I$.  Bounds for $b(I)$ and $c(I)$ are given in \cite{CHT}. 
Note  that if the  ideal $I$ is generated in a single degree, say $s$, then   $a(I)=s$
and hence $\reg(I^k)-\reg(I^{k-1}=s$ for large $k$. We say that the regularity of the
powers of $I$  jumps at place $k$ if  $\reg(I^k)-\reg(I^{k-1}>s$. 

One of the most powerful tools in proving that an ideal has a linear
resolution is the following notion:  

\begin{Definition} 
 An ideal $I$ generated in a single degree is said to have linear
quotients if there exists  a system of minimal generators $f_1,\dots, f_s$ of $I$ such
that for every $k\leq s$  the colon ideal
$(f_1,\dots, f_{k-1}):f_k$ is generated by linear forms.
\end{Definition} 

One has (see \cite{CH}): 

\begin{Lemma} 
\begin{itemize} 
\item[(a)] If $I$ has linear quotients then $I$ has a linear resolution. 
\item[(b)] If $I$ is a monomial ideal, then the property of having linear
quotients with respect to its monomial generators is independent of the
characteristic of the  base filed. 
\end{itemize} 
\end{Lemma}

\section{Examples} 

In this section we present some (known and some new)   examples of ideals with a linear
resolution such that the square does   have non-linear syzygies.  

The first example of such an ideal was discovered  by Terai. It is an
ideal well-known  for having another pathology: it is a
square-free monomial ideal whose Betti numbers,  regularity and projective
dimension depend on the characteristic of the base field. 

\begin{Example}\label{terai}     Consider the ideal 
$$I=(abc, abd, ace, adf, aef, bcf, bde, bef, cde, cdf)$$
 of $K[a,b,c,d,e,f]$. In
characteristic $0$  one has
$\reg(I)=3$ and $\reg(I^2)=7$. The only non-linear syzygy for $I^2$ comes
at the very end of the resolution. In characteristic $2$ the ideal $I$
does not have a linear resolution.  
\end{Example} 

The second example is taken from \cite{CH}. It is monomial and  characteristic free. 
The ideal is defined by  $5$ monomials and  very likely  there are no  such examples
with less than $5$ generators.

\begin{Example}\label{ConcaHerzog}     Consider the ideal 
$$I=(a^2b, a^2c, ac^2, bc^2, acd)$$
of $K[a,b,c,d]$ . It is easy to check that $I$ has
linear quotients (with respect to the monomial generators in the given
order). It follows that $I$  has a linear resolution independently  of
$\chara K$. Furthermore
$I^2$ has a quadratic first-syzygy in characteristic $0$.  But the first
syzygies of a monomial ideals are independent of $\chara K$. So we may
conclude that $\reg(I)=3$ and
$\reg(I^2)>6$ for every base field $K$.   
\end{Example}

The third example is due to Sturmfels \cite{S}. It is monomial, square-free  and 
characteristic free. It is  defined by  $8$ square-free monomials and, according
to Sturmfels \cite{S},   there are no such   examples with less than $8$
generators.

\begin{Example}\label{Sturmfels}     Consider the ideal  
$$I=(def, cef, cdf, cde, bef, bcd, acf, ade)$$
of $K[a,b,c,d,e,f]$.  One checks 
that $I$ has linear quotients (with respect to the monomial generators in the given
order) and so it has a linear resolution independently  of $\chara K$. Furthermore
$I^2$ has a quadratic first-syzygy. One concludes that $\reg(I)=3$ and
$\reg(I^2)>6$ for every base field $K$.   
\end{Example} 

So far all the examples were monomial ideals generated in degree $3$. 
Can we  find examples generated in degree $2$?   In view of the main result of
\cite{HHZ},   we have to allow also non-monomial generators. 
One binomial generator is enough:

 \begin{Example}  Consider the ideal 
$$I=(a^2,  ab,  ac, ad, b^2, ae + bd, d^2)$$
 of $K[a,b,c,d,e]$.  
One   checks   that  $\reg(I)=2$ and $\reg(I^2)=5$ in characteristic
$0$. Very likely the same holds in any characteristic. The ideal has
linear quotients with respect to the  generators in the given order. 
\end{Example}

One may wonder whether there exists a prime ideal with this behavior. 
Surprisingly, one can find such an example already among the most beautiful   and
studied prime ideals, the generic determinantal ideals. 

\begin{Example} \label{prime}    Let $I$ be the ideal of $K[x_{ij} : 1\leq
i\leq j\leq 4]$ generated by the  the $3$-minors of the generic symmetric
matrix $(x_{ij})$. It is well-known that $I$ is  a prime ideal defining a
Cohen-Macaulay ring and that $I$ has a linear resolution. One  checks  that $I^2$
does not have a linear resolution.  
\end{Example} 

\begin{Remark} {\rm
Denote by $I$ the ideal of \ref{prime} and by $J$ the ideal of
\ref{terai}.  It is interesting to note that the graded  Betti numbers 
of $I$ and $J$ as well as those  of $I^2$ and $J^2$ coincide.   Is this just an
accident? There might be some hidden relationship between the two ideals,
e.g.  $J$  could be an initial ideal or  a specialization (or an  initial
ideal of a specialization) of the ideal $I$. Concretely, we may ask
whether $J$ can be represented as the (initial) ideal of (the ideal of) 
$3$-minors of a
$4\times 4$ symmetric matrix of linear forms in $6$ variables. 
We have not been able to answer this question but we believe that
something like that should be true. Note however that the most natural way of
filling a $4\times 4$ symmetric matrix with $6$ variables would be to put $0$'s
on the main diagonal and to fill the remaining positions with the $6$
variables. Taking $3$-minors one gets an ideal, say $J_1$,  which shares
many invariants with $J$. For instance, we have checked that $J_1$ and
$J$ as well as their squares have the same graded Betti
numbers (respectively of course).  The ideals $J$ and $J_1$  are both reduced  but $J$
has $10$ components of degree $1$  while $J_1$ has  $4$ components of  degree $1$
and  $3$ of degree $2$. We have also checked that, in the given
coordinates, $J$ cannot be an initial ideal of $J_1$.   }
\end{Remark} 

\section{Regularity Jumps} 

The goal of this section is to   show that the regularity of the powers of  an ideal can
jump  for the first time at any place.  This happens already in $4$
variables.  Let $R$ be  the polynomial ring
$K[x_1,x_2,z_1,z_2]$  and for $d>1$ define the ideal 

$$J=(x_1z_1^d, x_1z_2^d, x_2z_1^{d-1}z_2)+(z_1, z_2)^{d+1}$$ 

We will prove the following: 

\begin{Theorem}
\label{main}  The ideal $J^k$ has linear quotients for all  $k<d$ and
$J^d$ has a first-syzygy of degree $d$. In particular 
$\reg(J^k)=k(d+1)$ for all $k<d$ and $\reg(J^d)\geq d(d+1)+d-1$ and this holds
independently of $K$. 
\end{Theorem} 

In order to prove that $J^k$ has linear quotients we need the following 
technical construction. 

Given  three sets of variables $x=x_1,\dots,x_m$, $z=z_1,\dots,z_n$ and
$t=t_1,\dots,t_k$ we consider  monomials
$m_1,\dots,m_k$ of degree $d$ is the variables $z$.  Let $\phi$ be    a
map
$$\phi:\{t_1,\dots,t_k\}\to
\{x_1,\dots,x_m\}. $$
Extend its action to arbitrary monomials by
setting   
$$\phi(\prod t_i^{a_i})=\prod \phi( t_{i})^{a_i}.$$  Set
$W=(m_1,\dots,m_k)$.   Consider the bigraded  presentation
$$\Phi:K[z_1,\dots,z_n,t_1,\dots,t_k] \to  R(W)=K[z_1,\dots,z_n ,m_1s,\dots, m_ks]$$ 
of the  Rees algebra  of the ideal $W$ obtained by setting  $\Phi(z_i)=z_i$ and
$\Phi(t_i)=m_is$ ($s$ a new variable) and giving degree $(1,0)$ to the $z$'s and degree
$(0,1)$ to the $t$'s. We set 

$$H=\Ker \Phi.$$
and note that $H$ is a binomial ideal.

\begin{Definition}\label{pseudo}   We say that
$m_1\dots,m_k$ are pseudo-linear of order
$p$ with respect to $\phi$ if for every $1\leq b\leq a\leq p$ and for
every binomial
$MA-NB$  in $H$    with $M,N$ monomials in the
$z$ of degree $(a-b)(d+1)$ and $A,B$  monomials in the $t$ of degree $b$
such that
$\phi(A)>\phi(B)$ in the lex-order there exists an element of the form
$M_1t_i-N_1t_j$ in $H$ where $M_1,N_1$ are monomials in the $z$ such
that the following conditions are satisfied: 
\begin{itemize}
\item[(1)]  $N_1|N$,  
\item[(2)]  $t_i|A$,  $t_j|B$,
\item[(3)]  $\phi(t_i)|\phi(A)/\GCD(\phi(A),\phi(B))$,
\item[(4)]  $\phi(t_i)>\phi(t_j)$ in the lex-order. 
\end{itemize} 
\end{Definition} 

The important consequence is the following: 

\begin{Lemma}
\label{pl1} Assume that  $m_1\dots,m_k$ are monomials of degree $d$ in the
$z$ which are pseudo-linear of order $p$ with respect to $\phi$. Set
$$J=(z_1,\dots,z_n)^{d+1}+(\phi(t_1)m_1,\phi(t_2)m_2,\dots,
\phi(t_k)m_k).$$ Then $J^a$ has linear quotients for all $a=1,\dots,p$. 
\end{Lemma} 

\demo  Set $Z=(z_1,\dots,z_n)^{d+1}$ and
$I=(\phi(t_1)m_1,\phi(t_2)m_2,\dots,
\phi(t_k)m_k)$. Take $a$ with $1\leq a\leq p$ and order the generators of
$J^a$ according to the following decomposition:
$J^a=Z^a+Z^{a-1}I+\dots+Z^bI^{a-b}+\dots+I^a$. In the block
$Z^a$ we order the generators so that they have linear quotients; this is
easy since $Z^a$ is just a power of the $(z_1,\dots,z_n)$. In the the
block $Z^bI^{a-b}$ with $b<a$   we order the generators extending (in
anyway)  the lex-order in the $x$. We claim that, with this order, the
ideal $J^a$ has linear  quotients. Let us check this. As long as we deal
with elements of the block $Z^a$ there is nothing to check. So let us
take some monomial, say $u$ from the block $Z^bI^{a-b}$ with $b<a$ and
denote by $V$ the ideal of generated by the monomials which are earlier
in the list. We have to show that the colon ideal  $V:(u)$ is generated
by variables. Note that  $V:(u)$ contains
$(z_1,\dots,z_n)$ since $(z_1,\dots,z_n)u\subset Z^{b+1}I^{a-b-1}\subset
V$. Let $v$ be a generator of $V$. If $v$ comes from a block
$Z^cI^{a-c}$ with $c>b$ then we are done since 
$(v):(u)$ is contained in $(z_1,\dots,z_n)$ by degree reason. So we can
assume that also
$v$ comes from the block $Z^bI^{a-b}$.  Again, if the generator
$v/\GCD(v,u)$ of $(v):(u)$ involves the variables
$z$  we are done. So we are left with the case in which  $v/\GCD(v,u)$
does not involves the variables $z$. It is now the time to use the
assumption that the
$m_i$'s are pseudo-linear.  Say $u=Nm_{s_1}\phi(t_{s_1}) \cdots
m_{s_a}\phi(t_{s_a})$ and
$v=Mm_{r_1}\phi(t_{r_1}) \cdots m_{r_a}\phi(t_{r_a})$ with $M,N$ monomials
of degree
$(b-a)(d+1)$ in the $z$. Set $A=t_{r_1} \cdots t_{r_a}$ and $B=t_{s_1}
\cdots t_{s_a}$.  
 Since $v$ is earlier than $u$ in the generators of $J^a$ we have 
$\phi(A) >\phi(B)$ in the lex-order. Note also that $(v):(u)$ is
generated by $\phi(A)/\GCD(\phi(A),\phi(B))$.
 Now the fact that 
$v/\GCD(v,u)$ does not involves the variables $z$ is equivalent to say
that 
$MA-NB$ belongs to $H$. By assumption   there exists 
$L=M_1t_{i}-N_1t_{j}$  
 in $H$ such that  the conditions (1)--(4) of Definition \ref{pseudo} hold. Multiplying
$L$ with $(N/N_1)(B/t_j)$, we have  that
$M_1(N/N_1)t_i(B/t_j)-NB$ is in $H$ and by construction
$$v_1=M_1 (N/N_1) m_i\phi(t_i) m_{s_1}\phi(t_{s_1}) \cdots
m_{s_a}\phi(t_{s_a}) / m_{j}\phi(t_{j})$$
 is a monomial of the block $Z^bI^{a-b}$ which
is in $V$ by construction and such that 
$(v):(u)\subseteq (v_1):(u)=(\phi(t_i))$. 
This concludes the proof. 
\QED

Now we can prove: 

\begin{Lemma}\label{pl2}  For every integer  $d>1$ the monomials 
$$m_1=z_1^d,\ \  m_2=z_2^d,\ \  m_3=z_1^{d-1}z_2$$  are pseudo-linear of
order $(d-1)$ with respect to the map 
$$\phi:\{t_1,t_2,t_3\}\to \{x_1,x_2\}$$ defined by $\phi(t_1)=x_1$,
$\phi(t_2)=x_1$,
$\phi(t_3)=x_2$. 
\end{Lemma} 

\demo It is easy to see that the defining ideal $H$  of  the Rees algebra of
$W=(m_1,m_2,m_3)$ is generated by  
$$(3)\  z_2t_1-z_1t_3\qquad (4)\  z_1^{d-1}t_2-z_2^{d-1}t_3 \qquad (5)\ 
t_1^{d-1}t_2-t_3^d.$$

Let $1\leq b\leq a\leq d-1$ and $F=MA-NB$ a binomial of bidegree
$((a-b)(d+1), b)$ in $H$ such that $\phi(A)>\phi(B)$ in the lex-order. 
Denote by
$v=(v_1,v_2), u=(u_1,u_2)$ the exponents of $M$ and $N$ and by 
$\alpha=(\alpha_1,\alpha_2,\alpha_3)$ and
$\beta=(\beta_1,\beta_2,\beta_3)$ the exponents of $A$ and $B$.  We
collect all the  relations that hold by assumption: 

$$\begin{array}{ll} (i)   & 1\leq b\leq a\leq d-1\\  (ii)  &
\alpha_1+\alpha_2+\alpha_3=\beta_1+\beta_2+\beta_3=b\\ (iii) & 
v_1+v_2=u_1+u_2=(a-b)(d+1)\\ (iv)  & 
v_1+d\alpha_1+(d-1)\alpha_3=u_1+d\beta_1+(d-1)\beta_3\\ (v)   & 
v_2+d\alpha_2+\alpha_3=u_2+d\beta_2+\beta_3\\ (vi)  & 
\alpha_1+\alpha_2>\beta_1+\beta_2\\
\end{array}
$$

Note that $(vi)$ holds since  $\phi(A)>\phi(B)$ in the lex-order. If  
$$   \alpha_1>0 \mbox{ and } \beta_3>0 \mbox{ and }  u_1>0 \eqno{(6)}$$
then equation (3) does the job.  If instead   
$$   \alpha_2>0 \mbox{ and } \beta_3>0  \mbox{ and }  u_2\geq d-1
\eqno{(7)}$$  then equation (4) does the job.  So it is enough to show
that either $(6)$ or $(7)$ hold. By contradiction, assume that both $(6)$
and $(7)$ do not hold. Note that $(ii)$ and $(vi)$ imply that
$\beta_3>0$. Also note that if $a=b$ then the equation $F$ has bidegree
$(0,a)$ and hence must be divisible by $(3)$ which is impossible since
$a<d$. So we may assume that
$b<a$. But then 
$u_1+u_2=(a-b)(d+1)\geq d+1$ and hence either $u_1>0$ or $u_2\geq d-1$.
Summing up, if both
$(6)$ and $(7)$  do not hold and taking into consideration that
$\beta_3>0$, that either 
$u_1>0$ or $u_2\geq d-1$ and that $\alpha_1+\alpha_2>0$, then one of the
following conditions hold: 

$$\begin{array}{ll} (8) &  \alpha_1=0 \mbox{  and  } u_2<d-1  \\ (9) & 
\alpha_2=0 \mbox{  and  } u_1=0 
\end{array}
$$

If $(8)$  holds then $\alpha_2>\beta_1+\beta_2$. Using $(v)$ we may write
$$   
u_2=v_2+d(\alpha_2-\beta_1-\beta_2)+d\beta_1+\alpha_3-\beta_3.\eqno{(10)}$$
If $\beta_1>0$ we conclude that $u_2\geq d+d-\beta_3$ and hence $u_2\geq
d+1$ since
$\beta_3\leq b\leq d-1$, and this is a contradiction. 

If instead $\beta_1=0$ then $\alpha_2+\alpha_3=\beta_2+\beta_3$ and hence
$(10)$ yields 
$u_2=(d-1)(\alpha_2-\beta_2)+v_2\geq (d-1)$, a contradiction.

If $(9)$  holds then $\alpha_1>\beta_1+\beta_2$. By $(iv)$ we have 
$$v_1+d(\alpha_1-\beta_1)+(d-1)(\alpha_3-\beta_3)=0$$ Hence
$$v_1+d(\alpha_1-\beta_1-\beta_2)+(d-1)(\alpha_3-\beta_3)+d\beta_2=0$$
But $\alpha_3-\beta_3=-\alpha_1+\beta_1+\beta_2$, so that 
$$v_1+(\alpha_1-\beta_1-\beta_2)+d\beta_2=0$$ which is impossible since
$\alpha_1-\beta_1-\beta_2>0$ by assumption. 
\QED

Now we are ready to complete the proof of the theorem. 
\medskip

\demo[ of Theorem \ref{main}]: Combining \ref{pl1} and \ref{pl2} we have
that
$J^k$ has linear quotients, and hence a linear resolution, for all
$k<d$.   It remains to show that $J^d$   has  a  first-syzygy of degree
$d$.
   Denote by
$V$ the ideal generated by all the monomial generators of $J^d$ but
$u=(z_1^{d-1}z_2x_2)^d$. We claim that
$x_1^d$ is a minimal generator of $V:u$ and this is clearly enough to
conclude that
$J^d$   has  a first-syzygy of degree $d$. First note that
$x_1^du=x_2^d(z_1^dx_1)^{d-1}(z_2^dx_1)\in  V$, hence $x_1^d\in V:u$.
Suppose, by contradiction, $x_1^d$ is not a minimal generator of $V:u$.
Then there exists an integer
$s<d$ such that $x_1^su\in V$. In other words we may write $x_1^su$ as
the product of $d$ generators of $J$, say $f_1,\dots,f_d$, not all equal
to $z_1^{d-1}z_2x_2$, times a monomial
$m$ of degree $s$.  Since the total degree in the $x$ in $x_1^su$ is
$s+d$ at each generator of $J$ has degree at most $1$ in the $x$, it
follows that   the $f_i$ are all of the of type $z_1^dx_1, z_2^dx_2,
z_1^{d-1}z_2x_2$ and $m$ involves only $x$.  Since
$x_2$ has degree $d$ in $u$,  $s<d$ and $z_1^{d-1}z_2x_2$ is the only
generator of $J$ containing $x_2$ it follows that at least one  of the
$f_i$ is equal  to $z_1^{d-1}z_2x_2$. Getting rid of those common factors
we obtain a relation of type
$x_1^s(z_1^{d-1}z_2x_2)^r=m(z_1^{d}x_1)^{r_1}(z_2^{d}x_1)^{r_2}$ with
$r=r_1+r_2<d$.  In the $z$-variables it gives
$ (z_1^{d-1}z_2)^r=(z_1^{d})^{r_1}(z_2^{d})^{r_2}$ with $r=r_1+r_2<d$
which is clearly impossible. 
\QED 

What is the regularity of $J^k$ for $k\geq d$?  There is some  computational evidence
that  the first guess, i.e.  $\reg(J^k)=k(d+1)+d-1$ for  $k\geq d$,  
might be correct . 
 
\section{Variations}

The ideas and the strategy of the previous section can be used, in
principle,  to create other kinds of ``bad"  behaviors. We give in this section
some hints and examples but no  detailed proofs.

\begin{Hint}  {\rm  Given $d>1$  consider the ideal 
$$H=(x_1z_1^d, x_1z_2^d, x_2z_1^{d-1}z_2)+z_1z_2(z_1, z_2)^{d-1}$$ 
We believe that $\reg(H^k)=k(d+1)$ for all $k<d$ and $\reg(H^d)\geq
d(d+1)+d-1$.  Note that $H$ has two generators less than
$J$. In the case $d=2$, $H$ is exactly the ideal of \ref{ConcaHerzog}. 
}
\end{Hint}

One can ask whether there are radical   ideals with  a behavior
as the ideal in  \ref{main}. One would need a square-free version of   the construction
of the previous section. This suggests the following:

\begin{Hint}  {\rm  For every $d$ consider variables $z_1,z_2,\dots,z_{2d}$
and $x_1,x_2$ and the ideal $$J=(
\begin{array}{ll}  x_1z_1z_2, x_1z_3z_4, \dots, x_1z_{2d-1}z_{2d},\\
x_2z_2z_3, x_2z_4z_5, \dots, x_2z_{2d}z_{1}
\end{array} )+\Sq^3(z)$$
 where $\Sq^3(z)$ denote the square-free cube of
$(z_1,\dots,z_{2d})$, i.e. the ideal generated by the  square-free monomials of degree
$3$ in the $z$'s. We conjecture that $\reg(J^k)=3k$ for $k<d$ and $\reg(J^d)>3d$.
 Note that for $d=2$ one obtains Sturmfels' Example \ref{Sturmfels}. }
\end{Hint}

We have no idea on how to construct  prime ideals with  a behavior as the ideal in
\ref{main}.   If one  wants two (or more) jumps one can try with:

\begin{Hint} {\rm Let $1<a<b$ be integers. Define  the ideal 
$$I=(y_2z_1^b, y_2z_2^b, xz_1^{b-1}z_2)+
z_1^{b-a} (y_1z_1^a, y_1z_2^a,xz_1^{a-1}z_2)+
z_1z_2(z_1,z_2)^{b-1}$$
of the polynomial ring $K[x,y_1,y_2, z_1,z_2]$. 
We expect that $\reg(I)=b+1$ and 
$\reg(I^k)-\reg(I^{k-1})>(b+1)$   if  $k=a$ or $k=b$.  
}\end{Hint}

\smallskip

\begin{center}

{\Large\bf Acknowledgments}

\end{center}
The author wishes to thank the organizers
of the  Lisbon Conference on Commutative
Algebra (Lisbon, June 2003)  for the kind invitation  and for their warm hospitality. 
Some parts of this  research project was carried out  while the author was visiting 
MSRI (Berkeley) within the  frame of  the  Special Program 2002/03 on Commutative
Algebra.   The results and examples    presented in  this  the paper  have been inspired
and suggested by  computations   performed  by   the computer algebra system CoCoA  
\cite{CNR}.

\end{document}